\documentclass[12pt]{article}
\usepackage{amsmath,latexsym,amsthm,amsfonts,amssymb, makeidx}
\date{}

\newtheorem{theorem}{Theorem}

\newtheorem{lemma}[theorem]{Lemma}

\newtheorem{proposition}[theorem]{Proposition}

\newtheorem{problem}[theorem]{Problem}

\newcommand{\IQ}{\mathbb Q}

\newcommand{\w}{\omega}

\newcommand{\U}{\mathcal{U}}
\newcommand{\Op}{\mathcal{O}\/}

\newcommand{\C}{\mathcal{C}}

\newcommand{\A}{\mathcal{A}}

\newcommand{\B}{\mathcal{B}}

\newcommand{\F}{\mathcal{F}}

\newcommand{\PP}{\mathcal P}

\newcommand{\shu}{\bigcup_{fin}(\mathcal O,\Gamma)}
\newcommand{\ssc}{\bigcup_{fin}(\mathcal O,\Omega)}
\newcommand{\snu}{\mathrm{S}_{fin}(\Gamma,\Omega)}

\newcommand{\sme}{\bigcup_{fin}(\mathcal O,\mathcal O)}

\newcommand{\Prr}{\mathsf P}
\newcommand{\G}{\mathsf G}

\title{Can a Borel group be generated by a Hurewicz subspace?}
\author{Lyubomyr Zdomskyy}
\begin{document}

\maketitle

\begin{abstract}
In this paper we  formulate three problems concerning topological properties
of sets generating Borel non-$\sigma$-compact groups.
In case of the concrete $F_{\sigma\delta}$-subgroup of
$\{0,1\}^{\w\times\w}$ this gives an equivalent reformulation
of the Scheepers diagram problem.
\end{abstract}

\noindent\large\textbf{Introduction}\normalsize
\medskip

The Hurewicz property
was introduced in \cite{Hu} as a cover counterpart
of the  $\sigma$-compactn\-ess: a topological space $X$ is said to have this 
property,
if for every sequence $(u_n)_{n\in\w}$ of open covers of $X$ there exists
a sequence $(v_n)_{n\in\w}$, where each $v_n$ is a finite subset of $u_n$,
such that each element $x\in X$ belongs to $\cup v_n$ for all  but finitely many
$n\in\w$. It is easy to see that each $\sigma$-compact space is Hurewicz (= has 
the Hurewicz property).
The converse statement is known to fail in ZFC, see \cite{JMSS}.
By a \emph{Borel} space we mean a separable metrizable space which is a Borel 
subset of its
completion.  This paper is devoted to problems close to the subsequent one.
\begin{problem} \label{pr1}
Can a Borel non-$\sigma$-compact group be generated by its  Hurewicz subspace?
\end{problem}
This problem is especially interesting for the concrete subgroup $\G$
of $\{0,1\}^{\w\times\w}$ (standardly endowed with the coordinatewise addition 
modulo $2$)
 being equivalent to the ``Hurewicz'' part of
the Scheepers diagram problem (see \cite[Problems~1,2]{JMSS}, 
\cite[Problems~4.1,4.2]{SPM2}, 
 \cite[Problem~1]{Ts-sur-p}, and \cite[Problem~3.2]{Ts-sur}), where
$$     \G  =\{x\in\{0,1\}^{\w^2}: \mbox{ for every }j\in\w\mbox{ and  for all 
but finitely many }i\in\w\ (  x_{i,j}=0)\}. $$
In order to formulate the Scheepers diagram problem we have to recall  some 
definitions.
M. Scheepers in his work \cite{Sc} introduced a long list of new properties
looking similar to the Hurewicz one, and thus gave rise to the branch of set-
theoretic
topology known as \emph{Selection Principles}.
Selection principles may be thought as some combinatorial conditions on the 
family
of  open covers of a topological space.
 Let $\A$ and $\B$ be a families of covers of a topological space $X$.
Following \cite{Sc} we say that  $X$ has the property
\begin{itemize}
\item $\bigcup_{fin}(\A,\B)$, if for every sequence $(u_n)_{n\in\w}\A^\w$ 
 there exists a sequence $(v_n)_{n\in\w}$, where each $v_n$ is a finite subset
of $u_n$, such that $\{\cup v_n:n\in\w\}\in\B$; 
\item $\mathrm{S}_{fin}(\A,\B)$, if for every sequence $(u_n)_{n\in\w}\in\A^\w$ 
there exists a sequence $(v_n)_{n\in\w}$, where each $v_n$ is a finite subset
of $u_n$, such that $\bigcup\{v_n:n\in\w\}\in\B$.
\end{itemize}
Throughout the paper $\A $ and $\B$ run over the families  $\Op$, $\Omega$, and
$\Gamma$ of all open ($\w$-, $\gamma$-) covers of $X$.
Given a family $u=\{U_i:i\in I\}$ of subsets of a set $X$, we define
the map $\mu_u:X\to\PP(I)$ letting $\mu_u(x)=\{i\in I:x\in U_i\}$
($\mu_u$ is nothing else but the  Marczewski ``dictionary'' map introduced in 
\cite{M-S}). 
In what follows $I\in\{\w,w^2\}$. 
Depending on the properties of $\mu_u(X)$ a family $u=\{U_n:n\in\w\}$ is defined 
to be
\begin{itemize}
\item an \emph{$\w$-cover} \cite{GN}, if the family $\mu_u(X)$ is centered,
i.e. for every finite subset $K$ of $X$ the intersection $\bigcap_{x\in 
K}\mu_u(x)$
is infinite; 
\item a \emph{$\gamma$-cover} of $X$ \cite{GN}, if for every $x\in X$ the set 
$\mu_u(x)$
is cofinite in $\w$, i.e.  $\w\setminus \mu_u(x)$ is finite.
\end{itemize}
We shall consider here  four selection principles:
$\shu$, $\ssc$, $\sme$ and $\snu$. Let us note that $\shu$ is nothing else but 
the
Hurewicz property. Concerning $\sme$, it is the classical Menger covering 
property
introduced in \cite{Me}. We are in a position now to formulate the
\medskip

\noindent\textbf{Scheepers diagram problem.} {\it \begin{itemize}
\item[$(1)$] Does the property
$\ssc$ imply $\snu$? 
\item[$(2)$] And if not, then does  $\shu$ imply $\snu$?
\end{itemize}}

\medskip

One may ask the same question as in Problem~\ref{pr1} for properties
$\ssc$ and $\sme$. 
\begin{problem} \label{pr2}
Can a Borel non-$\sigma$-compact group be generated by its subspace
with the property $\ssc$?
\end{problem}
\begin{problem} \label{pr3}
Can a Borel non-$\sigma$-compact group be generated by its subspace
with the property $\sme$?
\end{problem}
The subsequent theorem, which is the main result of this paper, is the 
reformulation
of a Scheepers diagram problem in algebraic manner.

\begin{theorem} \label{tmain}
The property $\shu$ (resp. $\ssc$, $\sme$) implies $\snu$ if and only if the group
$\G$ is not generated by its subspace with the property $\shu$ (resp. $\ssc$, $\sme$). 

In other words, the positive answer onto the Scheepers diagram problem (1) 
(resp. (2))
is equivalent to the negative answer onto Problem~\ref{pr2} (resp. 
Problem~\ref{pr1})
in case of the group~$\G$.  
\end{theorem}
The group $\G$ is a rather simple object from  the point of view
of Descriptive Set Theory. For every $j\in\w$  its projection onto 
$\{0,1\}^{\w\times\{j\}}$
is homeomorphic to $\IQ$ being a countable metrizable space without isolated 
points.
From the above it follows that  $\G$ is a countable intersection of $F_\sigma$
subsets of $\{0,1\}^{\w^2}$ (i.e. it is an $F_{\sigma\delta}$- or, equivalently, 
$\Pi^0_3$- subset) homeomorphic
to $\IQ^\w$. Therefore it is a nowhere locally-compact, and it fails to have the 
property $\sme$.
 For more simple groups from the point of view of Borel hierarchy
 Problem~\ref{pr1} can be answered in negative.
\begin{proposition} \label{1negative}
 No Borel non-$\sigma$-compact group $B$ can  be generated by its  subspace
$X$  with the property $\shu$ provided $B$ is an $F_\sigma$- or $G_\delta$-
subspace of a complete
metric space.
\end{proposition}
Recall that  a map $f$ from a topological space $X$ to a topological space $Y$
is \emph{Borel,} if for every Borel subset $B$ of $Y$ its preimage $f^{-1}(B)$
is a Borel subset of $X$. The subsequent statement answers Problem~\ref{pr3}
in positive under the Continuum Hypothesis. On the other hand, it is known that 
the properties $\ssc$ and $\sme$ coincide in some models of ZFC, see \cite{Zd}.
Therefore the negative answer onto Problem~\ref{pr2} would imply that the negative 
answer onto Problem~\ref{pr3} is consistent as well.
\begin{proposition}\label{3positive}
Under the Continuum Hypothesis a metrizable separable  group $B$ can be 
generated by its subspace
$X$ with the property $\sme$ provided it is a Borel homomorphic image of a 
 nonmeager metrizable separable group.  In particular, $\G$ is generated by its 
subspace
with the property $\sme$ under CH.
\end{proposition}  

\noindent\textbf{Remark.}
None of the known methods of contruction of spaces
with the property $\shu$ can give a subspace of a Borel non-$\sigma$-compact group
generating it. All finite
powers of  spaces with the property $\shu$ 
constructed in \cite[Theorem~5.1]{JMSS}, \cite[Theorem~5.1]{TZ}, and 
\cite[Theorem~10(1)]{BT}
have the property $\sme$ or even $\shu$, and hence so is any group they 
generate.  
But every Borel (even analytic) space with the property $\sme$ is $\sigma$-compact,
see \cite{Ar}.
While the Sierpinski sets $S$ considered in \cite{JMSS} and \cite{ST} have the 
subsequent
 property: for every Borel subset $B$ containing $S$ there exists a $\sigma$-compact $L$
such that $S\subset L\subset B$, see \cite{BZ}.

Concerning the property $\ssc$, all known examples (besides Sierpinski sets)
have the property $\sme$ in all finite powers, and hence can not generate
non-$\sigma$-compact Borel group.
\hfill $\Box$
\medskip

\ 
\medskip

\noindent\large \textbf{Proofs} \normalsize
\medskip

In  what follows $A\subset^\ast B$
standardly means that $A\setminus B$ is finite.
In our proofs we shall  exploit set-valued maps. By a \emph{set-valued map} $\Phi$ from a 
set $X$
into a set $Y$ we understand a map from $X$ into $\mathcal{P}(Y)$ and write
$\Phi:X\Rightarrow Y$ (here $\mathcal P(Y)$ denotes the set of all subsets of 
$Y$). For a subset $A$ of $X$ we put $\Phi(A)=\bigcup_{x\in A}\Phi(x)\subset Y.$
The set-valued map $\Phi$ between topologival spaces $X$ and $Y$ is said to be
\begin{itemize}
\item \emph{compact-valued}, if $\Phi(x)$ is compact for every $x\in X$;
\item\emph{upper semicontinuous}, if for every open subset $V$ of $Y$ the set
$\Phi_\subset^{-1}(V)=\{x\in X:\Phi(x)\subset V\}$ is open in $X$.
\end{itemize}
For a set $X$ we can identify $\PP(X)$ with the compact space $\{0,1\}^X$
via the map $X\supset A\mapsto\chi_A\in\{0,1\}^X$ assigning to a subset
of $X$ its characteristic function.
A family $\A$ of subsets of a set $X$ is called \emph{upward closed},
for every $A\in \A$ and $B\supset A$ we have $B\in\A$.
For a set $A\subset X$ we make the subsequent notation: $\uparrow A=\{B\subset 
X:A\subset B\}$.
The following lemma is a more convenient reformulation of Theorem~\ref{tmain}.
\begin{lemma} \label{mes2} Let $\Prr$ be a topological property preserved by 
images under
upper semicontinuous compact-valued maps. Then
the following conditions are equivalent:
\begin{itemize}
\item[$(1)$] The  property $\Prr$  implies $\snu$;
\item[$(2)$] for every (upward-closed)  $\mathcal F\subset \PP(\w^2)$ with the 
property $\Prr$
 such that  $\omega\times \{j\}\subset^\ast F$ for every $F\in\mathcal F$
and $j\in\w$,  there exists a
sequence $(K_j)_{j\in\omega}$ of finite subsets of $\omega$ such that each 
element
of the smallest filter containing $\mathcal F$ meets 
$\bigcup_{n\in\omega}K_j\times\{j\}$.
\end{itemize}
\end{lemma}
\begin{proof} $(1)\Rightarrow(2)$. It simply follows from definition of the 
property
$S_{fin}(\Gamma,\Omega)$ and the observation that $\{\{F\in\F:F\ni 
(i,j)\}:i\in\w\}$
is an open  $\gamma$-cover of  $\F$ for every $j\in\w$.

$(2)\Rightarrow(1)$. Let $X$ be a topological space with the property $\Prr$ 
 and $(u_j)_{j\in\w}$
be a sequence of open $\gamma$-covers of $X$. Let us write $u_j$ in the form
$u_j=\{U_{i,j}:i\in\w\}$. Set $u=\{U_{i,j}:i,j\in\w\}$.
Consider the set-valued map $\Phi:X\Rightarrow\PP(\w^2)$,
$\Phi:x\mapsto\uparrow\mu_u(x)$. Applying  Lemma~2 of \cite{Zd}, we conclude 
that
$\Phi$ is compact-valued and upper semicontinuous,
and hence $\F:=\Phi(X)$ has the property $\Prr$. The   definition of $\Phi$ 
implies that $\F$
is upward closed. Since $u_j$ is a $\gamma$-cover of $X$ for every $j\in\w$,
$\w\times\{j\}\subset^\ast F$ for each $F\in\F$. 
 From the above it follows that there exists a sequence $(K_j)_{j\in\w}$
of finite subsets of $\w$ such that each element of the smallest filter $\U$
containing $\F$ meets some $K_j\times\{j\}$. Then
the family $\{U_{i,j}:i\in K_j\}$ is easily seen to be an $\w$-cover of $X$, 
which finishes our proof.
\end{proof}   

The  properties $\sme$, $\ssc$,
and $\shu$ satisfy the conditions of the above lemma by \cite[Lemma~1]{Zd}.
\medskip

\noindent\textbf{Proof of Theorem~\ref{tmain}.} Let $\Prr$ be any of
the properties $\sme$, $\ssc$, and $\shu$. Assuming that $\Prr$ implies
$\snu$, fix a subspace $X$ of $\G$ with the property $\Prr$.
Let us denote by $\varphi$ the map assigning to a subset $A$ of $\w^2$
its characteristic function $\chi_A\in\{0,1\}^{\w^2}$.
Then the space $\F=\{\w^2\setminus A:A\in\varphi^{-1}(X)\}$
has  the property $\Prr$ being homeomorphic to $X$,
and $\w\times\{j\}\subset^\ast F$ for every $F\in\F$ by our choice of $\G\supset 
X$.
Applying Lemma~\ref{mes2}, we conclude that there exists a sequence 
$(K_j)_{j\in\w}$
of finite subsets of $\w$ such that $\bigcup_{j\in\w}K_j\times\{j\}$ meets all
elements of the smallest filter containing $\F$. Now, a direct verification 
shows that the charateristic function
$\chi_{\cup_{j\in\w}K_j\times\{j\}}$ can not be represented as a sum
of elements of $X$, which means that $X$ does not generate $\G.$

Next, let us assume that $\Prr $ does not imply $\snu$ and apply 
Lemma~\ref{mes2}
to find an upward closed family  $\F$ of subsets of $\w^2$
such that for every sequence $(K_j)_{j\in\w}$ of finite subsets of $\w$
there exists a finite subset $\A$ of $\F$ such that 
$$ (\bigcup_{j\in\w}K_j\times\{j\})\cap \bigcap \A=\emptyset. $$
Set $X=\{\chi_{\w^2\setminus F}:F\in\F\}$. Then $X$ has the property $\Prr$
being homeomorphic to $\F$. We claim that $X$ is a set of generators of $\G$.
Indeed, let us fix any $g\in \G$ and set $K_j=\{i\in\w:g_{i,j}=1\}$.
Then each $K_j$ is finite by the definition of $\G$. For the sequence 
$(K_j)_{j\in\w}$
find a finite subset $\A=\{A_i:i\leq n\}$ of $\F$ as above.
Using the upward closedness of $\F$,
define inductively a finite subset $\B=\{B_i:i\leq n\}$ of $\F$ 
letting $B_0=A_0$ and $B_k= A_k\cup\bigcup_{l<k} (\w^2\setminus B_l) $
for all $0<k\leq n$. It is easy to prove by induction over $k\leq n$
that  $(\w^2\setminus B_l)\cap(\w^2\setminus B_k)=\emptyset$ for all $l<k$
and $\bigcap_{l\leq k}B_k=\bigcap_{l\leq k}A_k$, 
consequently 
$\bigcap\B=\bigcap\A\subset(\w^2\setminus\bigcup_{j\in\w}K_j\times\{j\})$.
Let $C_k=B_k\cup(\w^2\setminus\bigcup_{j\in\w}K_j\times\{j\})$, $k\leq n$.
Then $\C=\{C_k:k\leq n\}$ has the following properties
\begin{itemize}
\item[$(i)$] $\cup\C=\w^2\setminus\bigcup_{j\in\w}K_j\times\{j\}$;
\item[$(ii)$] $(\w^2\setminus C)\cap (\w^2\setminus D)=\emptyset$ for all 
$C,D\in\C$;
\item[$(iii)$] $\C\subset\F$.
\end{itemize}
It sufficies to note that $\{\chi_{\w^2\setminus C_k}:k\leq n\}\subset X$ by 
$(iii)$
and $\chi_{\w^2\setminus C_0}+\cdots+\chi_{\w^2\setminus 
C_n}=\chi_{\cup_{j\in\w}K_j\times\{j\}}=g$,
which finishes our proof.
\hfill $\Box$
\medskip

\noindent\textbf{Proof of Proposition~\ref{1negative}.}
First assume that $B$ is a non-$\sigma$-compact $G_\delta$-subspace of a 
complete metric space
and fix a subspace $X$ of $B$ with the property $\shu$. The same argument as in
\cite[Theorem~5.7]{JMSS} gives  a $\sigma$-compact subset $L$ of $B$
such that $X\subset L$. Since $B$ is not $\sigma$-compact, it is not generated 
by
$L$, and hence by $X$ as well.

Now consider a non-$\sigma$-compact Borel group $B$ which is an $F_\sigma$-
subset
of a complete metric space $Y$ and write $B$ in the form $\bigcup_{n\in\w}B_n$,
where each $B_n$ is closed in $Y$. Let $X$ be a subspace of $B$ with the 
property
$\shu$. Since  the property $\shu$ is preserved by
closed subspaces, $X\cap B_n$ has the property $\shu$ for all  $n\in\w$.
In addition, each $B_n$ is a $G_\delta$-subspace of $Y$ being closed.
From the above it follows that there exists a $\sigma$-compact $L_n$
such that $X\cap B_n\subset L_n\subset B_n$, and consequently
$X\subset \bigcup_{n\in\w}L_n\subset B$. It sufficies to apply the
 same argument as in the first part of the proof.
\hfill $\Box$ 
\medskip

\noindent\textbf{Proof of Proposition~\ref{3positive}.}
Let $C$ be a nonmeager metrizable separable topological group  
and $f:C\to B$ be a surjective Borel homomorphism. Almost literal repetition of 
the proof
of Lemma~29 from \cite{ST} give us a subspace $Z$ of $ C$ such that $Z$ 
generates $C$
and each Borel image of $Z$ has the property $\sme$, see 
\cite[Corollary~30]{ST}.
It sufficies to note that $B$ is generated by $f(Z)$. 

Next, let us show that under CH the group $G$ is generated by its subspace with 
the property $\sme$.
Indeed, let us denote by  $\tau$  the Tychonoff product topology
on $\{0,1\}^{\w^2}=\prod_{j\in\w}\{0,1\}^{\w\times\{j\}}$, where  $\{0,1\}^{\w\times\{j\}}$ is considered with the 
discrete
topology for each $j\in\w$. Then $\tau|\G$ is stronger than
the natural topology on $\G$, and
 $(\G,\tau|\G)$ is a completely metrizable topological group being a countable
product of  countable discrete groups.
\hfill $\Box$
\medskip

\noindent \textbf{Acknowledgements.} The author wishes to express his thanks to 
professor Taras Banakh for supervising the writing of this paper.

\noindent Department of Mechanics and Mathematics,\\
 Ivan Franko Lviv National University, \\
Universytetska 1, Lviv, 79000, Ukraine.
\medskip

\textit{E-mail address:}   \texttt{lzdomsky@rambler.ru}
\end{document}